\documentclass[letterpaper]{article}
\usepackage{aaai}
\usepackage{times}
\usepackage{algorithm}
\usepackage{algpseudocode}
\usepackage{helvet}
\usepackage{courier}
\usepackage{float} 
\usepackage{amsmath,amssymb}
\usepackage{caption}
\usepackage{graphicx}
\usepackage{multirow}
\usepackage{breqn}
\usepackage{hyperref}
\usepackage{cite}
\usepackage[numbers]{natbib}
\usepackage{amsmath}

\hypersetup{
    colorlinks=true,
    linkcolor=blue,
    filecolor=magenta,      
    urlcolor=cyan,
    pdftitle={Lithium-Ion Battery Charging Schedule Optimization to Balance Battery Usage and Degradation},
    pdfpagemode=FullScreen}

\begin{document}

\title{Lithium-Ion Battery Charging Schedule Optimization to Balance \\ Battery Usage and Degradation}
\author{Jacob Azoulay\ and Nico Carballal\\
Stanford University\\
AA222: Engineering Design Optimization\\
jazoulay@stanford.edu  |  nicocarb@stanford.edu\\}
\maketitle

\begin{abstract}
\begin{quote}
This work optimizes a lithium-ion battery charging schedule while considering a joint revenue and battery degradation model. The study extends the work of Meheswari et. al. to encourage battery usage/charging at optimal intervals depending on energy cost forecasts. This paper utilizes central difference Nesterov momentum gradient descent to come to optimal charging strategies and deal with the non-linearities of the battery degradation model. This optimization strategy is tested against constant, random varied price forecasts and a novel Gaussian process cost forecasting model. Contrary to many other papers regarding battery charging, formulating schedule optimization as a multivariate optimization problem provides meaningful insight to the inherent balance between these two competing objectives.

\end{quote}
\end{abstract}

\section{Problem Statement} 
With increased electrification of devices and machinery used in our daily lives, there is an increased importance placed on energy storage, usage efficiency, and replenishment. Batteries are ubiquitous, varying in size, capacity, and type. In an ideal world, appliances that run on batteries would be able to transfer 100\% of their energy to power the device. Furthermore, the batteries would hold a maximal charge with no leakage overtime and the maximum charge capacity would not degrade over the lifetime of the battery. An ideal battery would also be capable of instantly recharging.

Of course, batteries are prone to inefficiencies; however, it is possible to optimize battery life, performance, and recharging time within certain constraints. The primary optimization objective considered in this paper is battery degradation, which is influenced primarily by charge rate and state of charge (SOC). Furthermore, to formulate a comprehensive model which encourages battery use, a competing energy cost minimization objective will be introduced, transforming the problem into a multi-objective optimization problem.

Various related works have been conducted analyzing the relationship between several battery performance specifications. Lei et al., for instance, determined a tradeoff between lithium-ion battery degradation rate and energy loss using first order approximations of energy dissipation \cite{lei_charging_2018}. Similarly, Maheshwari et al. determine a battery degradation model that can be used to predict battery performance depending on specified battery parameters \cite{maheshwari_optimizing_2020}. Several other papers analyze the correlation between battery performance and battery material, usage environment, and various other factors \cite{fayaz_optimization_2022} \cite{masias_opportunities_2021}. Even factors such as weather \cite{gharehghani_effect_2022} and general maintenance \cite{noauthor_factors_2013} have been examined in context of battery performance.

\section{Motivation}

While much research has been conducted in maximizing battery efficiency, very few papers have developed a robust model of short or long term degradation of lithium-ion batteries, let alone consider battery degradation in the first place. The study performed in this paper is greatly inspired by the work of Maheshwari et al. which introduced a lithium-ion battery degradation model based on empirical data and testing. Our paper, however, attempts to solve the proposed problem using different optimization approaches which are robust to non-linearities in the model \cite{maheshwari_optimizing_2020}.

Long term battery capacity degradation for lithium-ion batteries can be attributed to exogenous factors, such as environment temperature and age, and endogenous factors, such as SOC and operating current. Because oftentimes external factors contributing to battery degradation are more difficult to control for, this paper only considers internal factors when modeling the degradation. 

\begin{figure} [htbp]
    \centering
    \includegraphics[width=0.25\textwidth]{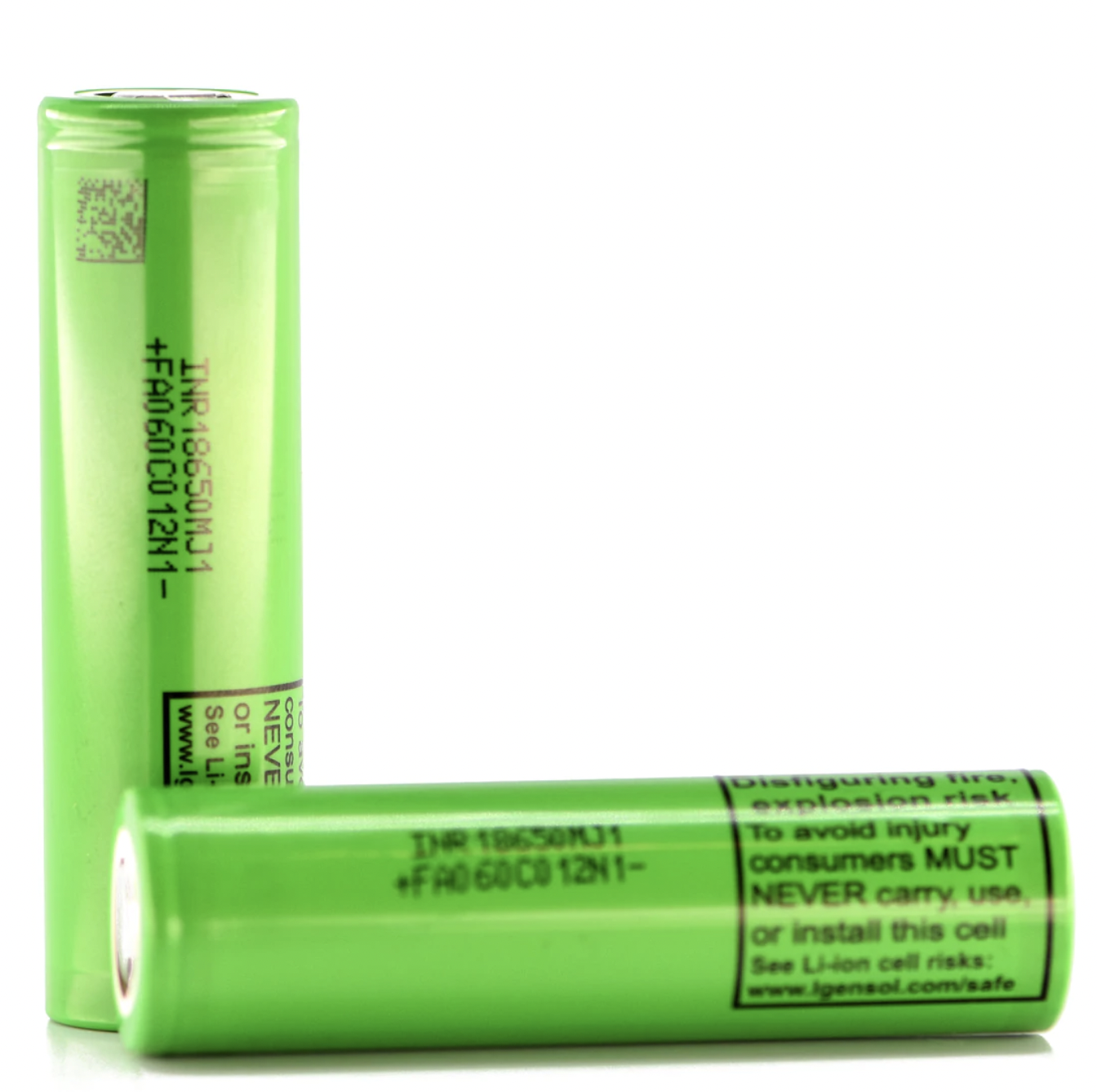}
    \caption{Two lithium-ion Sony US18650V3 batteries examined in this study \cite{maheshwari_optimizing_2020}.}
    \label{fig:battery_pic}
\end{figure}

Maheshwari et al. use experimental data of Sony US18650V3 lithium-ion batteries, pictured in Figure \ref{fig:battery_pic}, to demonstrate the correlation between capacity degradation and charge rate \cite{maheshwari_optimizing_2020}. Additionally, the paper showed the interwoven non-linear dependence between SOC, charge rate, and degradation. Because batteries are typically retired when their maximum capacity decreases to 80\%, better understanding these relationships will help model and optimize for prolonging battery lifetimes.

\section{Approach}
The goal of our project is to determine a battery charging/usage schedule that minimizes the long-term battery capacity degradation and the overall cost of charging the battery while still using the battery. As such, the problem will be modeled as a multi-objective non-linear optimization problem. To formulate a charging schedule, a finite time horizon between one-half day to two days, discretized into one-hour time steps, will be considered. At each time step, the design variables are charge (or discharge) rate $P_t^b$ for that given time interval. The final model output will be a sequence of charge rates which, if implemented at each corresponding time step, will minimize battery capacity degradation and minimize the cost of charging the battery according to the objective function.

The objective function of interest will be formulated as follows:
\begin{equation}
f(\textbf{P}) = w_C\times C + (1-w_C)\times D
\label{eq:objective-fn}
\end{equation}

where $C$ is total cost specified in the Cost Model Section, $D$ is total degradation specified in the Battery Degradation Model section, $w_C$ is the cost weighting factor, and $1 - w_C$ is the degradation weighting factor. The goal is to minimize both cost and degradation. A Pareto frontier can then be produced by sweeping through various values of $w_C \in [0, 1]$.

\subsection{Battery Degradation Model}

Because the charging schedule is composed of discrete time steps, total degradation $D$ will be the sum of the degradation experienced at each time step:

$$D = \sum_{t}^{T} d_t$$

where $d_t$ represents the degradation at time $t \in \{1, ..., T\}$.

The degradation model described below is inspired by \cite{maheshwari_optimizing_2020}. Maheshwari et. al have relevant tables and a much more in-depth description of each variable. This paper matches the conventions used by \cite{maheshwari_optimizing_2020}.

The SOC of the battery is updated at each time interval depending on the SOC of the previous state and on the total charge.

$$SOC_t = SOC_{t-1} + \frac{P_t^b \times \Delta T}{V_{nom} \times I_{1C}}$$

$P_t^b$ is the battery charge rate at time $t$, $\Delta T$ is the time interval of each step, $V_{nom}$ is the nominal operating voltage assumed to be constant, and $I_{1C}$ is the battery's 1C rate also assumed to be constant. Values for $\Delta T$, $V_{nom}$, and $I_{1C}$ used in this study are tabulated in Table \ref{tab:batParms}. 

\begin{table} [htbp]
\centering
\begin{tabular}{ |p{3cm}||p{3cm}|  }
 \multicolumn{2}{}{} \\
 \hline
 Parameter & Value\\
 \hline
 $\Delta T$    &  1 hr  \\
 $V_{nom}$     &  3.7 V \\
 $I_{1C}$      &  2.15 A\\
 $P^{b, max}$ &  5 W\\
 \hline
\end{tabular}
 \caption{Model parameters assumed in this study.}
 \label{tab:batParms}
\end{table}

To preserve battery health, manufacturers often recommend charging below a maximum charging rate $P^{b, max}$. Additionally, by definition, a battery's SOC must lie between zero and one-hundred percent. As such, $P_t^b$ and SOC will be constrained such that 
\begin{itemize}
    \item$\lvert P_t^b \rvert \leq P^{b, max} = 5 W$ 
    \item$0 \leq SOC_t \leq 100$
\end{itemize}

for all $t$.

Degradation for the battery considered in this study is highly dependant on SOC at the time of charge. For instance, charging a battery from $80\%$ to $100\%$ will result in greater degradation than charging the battery from $40\%$ to $60\%$ assuming the charge rate is constant. Because the degradation scheme is highly non-linear and battery specific, the degradation at varying SOC has been determined experimentally by Maheshwari et al. \cite{maheshwari_optimizing_2020}. The values recorded assume a constant current rate of 1C. The degradation at a 1C rate $d_t^{1C}$ at intermediate SOC can thus be interpolated from the recorded data.

\begin{figure} [htbp]
    \centering
    \includegraphics[width=0.45\textwidth]{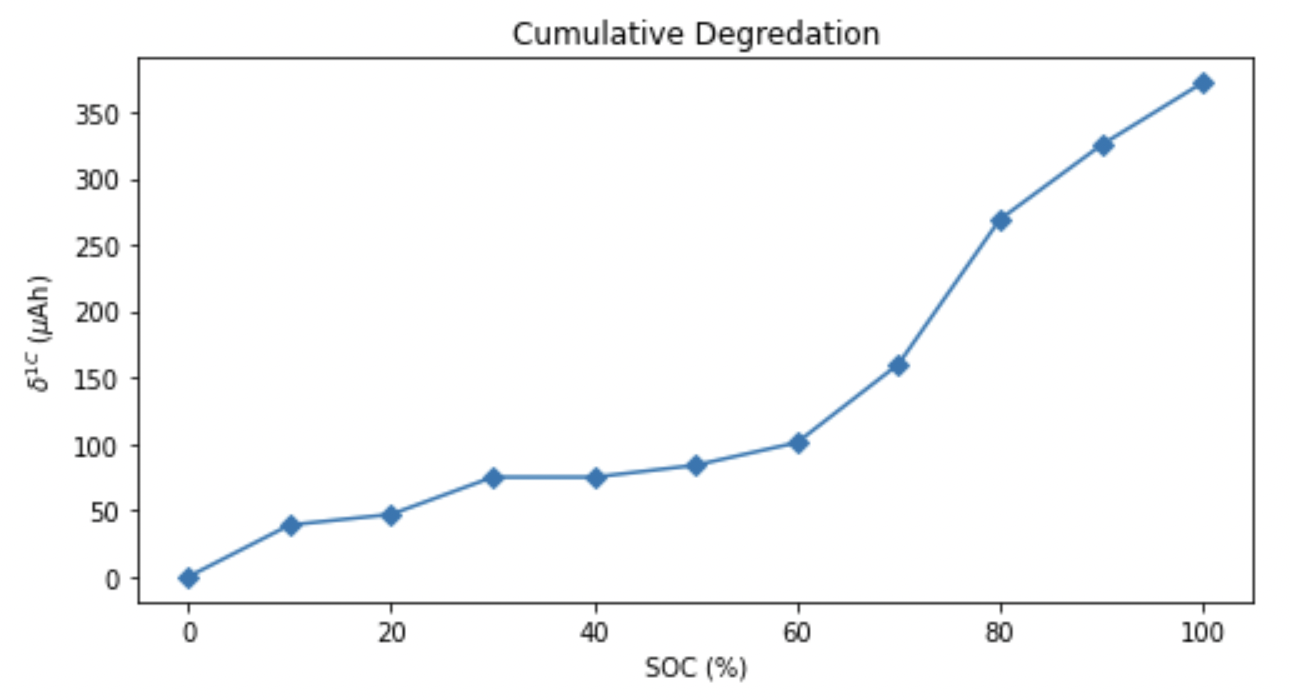}
    \caption{Cumulative degradation function determined experimentally to calculate $d_t^{1C}$.}
    \label{fig:cum_deg}
\end{figure}

To facilitate the evaluation of $d_t^{1C}$, the authors of \cite{maheshwari_optimizing_2020} developed a cumulative distribution graph from which $d_t^{1C}$ can easily be calculated. The graph in Figure \ref{fig:cum_deg} is a function of SOC on the x-axis and cumulative degradation on the y-axis. The equation to extract $d_t^{1C}$ is:

$$d_t^{1C} = \lvert \delta_t^{1C} - \delta_{t-1}^{1C} \rvert$$

\begin{figure} [htbp]
    \centering
    \includegraphics[width=0.45\textwidth]{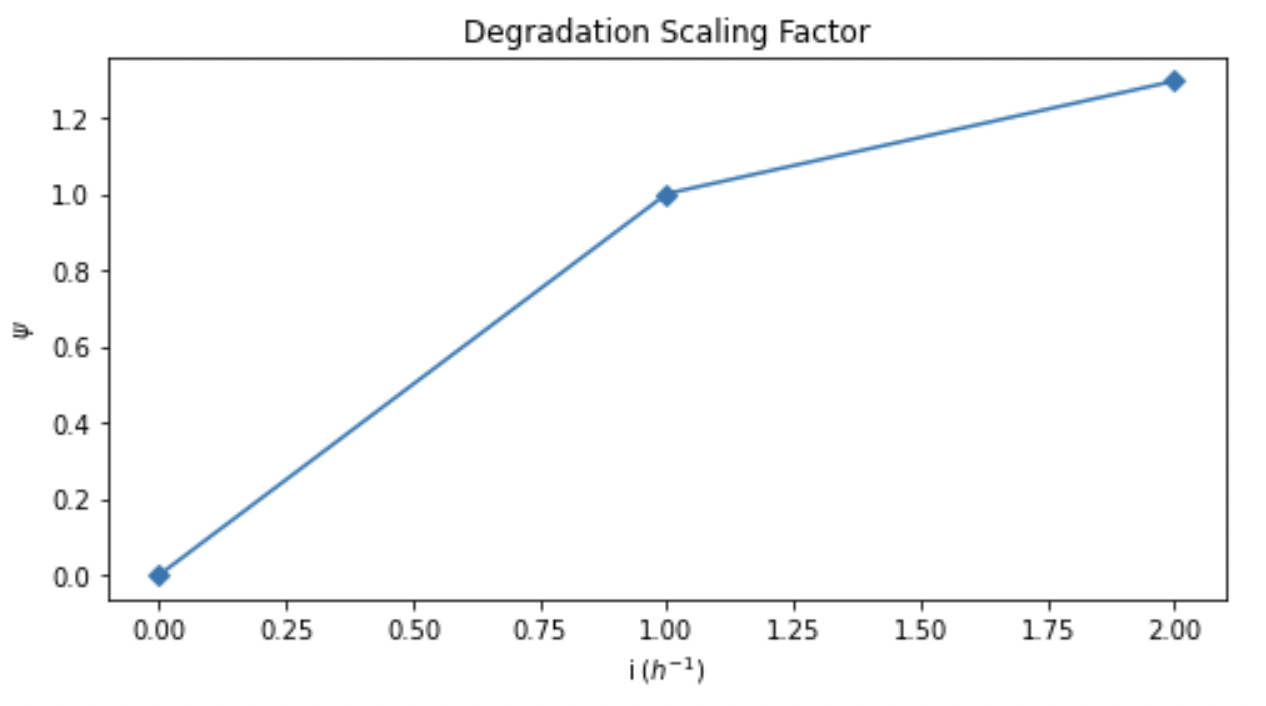}
    \caption{Current dependent scaling factor $\psi$ determined experimentally. Intermediary values for $\psi$ are calculated using linear interpolation. }
    \label{fig:psi_func}
\end{figure}

The obtained $d_t^{1C}$ value must then be scaled depending on current $i_t$. Higher current leads to higher degradation, therefore a scaling factor $\psi_t$ correlated to current is introduced, accounting for this dependence. At $i_t = 0$ (0C), $\psi_t=0$, indicating no degradation at zero load. At a current rate of 1C, $\psi_t= 1 $ by definition. Determined experimentally, $\psi_t=1.2956$ at a current rate of 2C. The model will estimate intermediary values of $\psi_t$ using linear interpolation. The graph of $\psi_t$ as a function of current rate can be seen in Figure \ref{fig:psi_func}. The mathematical formulations are as follows:

$$i_t = \frac{\lvert P_t^b \rvert}{V_{nom} \times I_{1C}}$$

$$
\psi_t(i_t)=
    \begin{cases}
        i_t & \text{if } 0 \leq i \leq 1\\
        1 + 0.2956 i_t & \text{if } 1 \leq i \leq 0
    \end{cases}
$$

Final degradation incurred $d_t$ is determined using the following equation:

$$d_t = \left(\frac{d_t^{1C} + \psi_t}{2}\right)^2 - \left(\frac{d_t^{1C} - \psi_t}{2}\right)^2$$

\subsection{Cost Model}

The cost model from the objective function is formulated to encourage battery usage
$$C = \sum_{t}^{T} P_t^b \times \lambda_t$$
where $P_t^b$ is the power to/from the battery during a trading interval. A trading interval, $t,$ is a period of time in which the price of energy, $\lambda_t$, is fixed. For this paper, trading intervals of an hour are considered --- meaning that prices are discretely changed at every hour interval. 

Because energy usage is inherently inefficient, an inefficiency factor, $\eta = .95$, is selected such that
$P_t^{b,ch} = \eta P_t^{m,ch}$ and $P_t^{b,dis} \eta = P_t^{m,dis}$, where $P_t^m$ is the power input/output of the market and $P_t^b$ is the power input/output directly of the battery. 

This cost model necessitates the price of energy over the entire trading period, $T$, at each trading interval, $\lambda_t$. A baseline constant price, $\lambda_t = p_{const}$ is tested first to ensure proper convergence. Next, random prices are tested. Lastly, a Gaussian Process fitting and sampling strategy is tested in the next section to predict and optimize to real data.

\subsection{Cost Model Prediction} \label{cost-pred}

The application of Gaussian processes for energy forecasting has been well-studied to allow for more advanced energy management techniques \cite{lubbe_evaluating_2020}\cite{lubbe_evaluating_nodate}. 

Given previous days' energy price data, this paper aims to model the next day's forecast for energy prediction while also quantifying uncertainty in the prediction \cite{mykel}. For a given current hour $h_i$, we accumulate the previous $m$ data points, $\boldsymbol{H} = [h_{i-m}, \dots, h_{i-1}]$ for hourly energy usage in MW, $\boldsymbol{E} = [E_{i-m}, \dots, E_{i-1}]$, from \cite{noauthor_hourly_nodate}. We aim to predict the next $n$ hours, $\boldsymbol{H^{*}} = [h_{i}, \dots, h_{i+n}]$,  energy usage to base our optimization on $\boldsymbol{\hat{E}} = [\hat{E_i}, \dots, \hat{E_{i+n}}]$ along with quantifying our uncertainty.

The posterior distribution of our model can be formulated as follows:
\begin{align} \label{eq:post}
    \boldsymbol{\hat{E}}\mid\boldsymbol{E} \sim
    \mathcal{N} (\boldsymbol{\mu}), \boldsymbol{\Sigma})
\end{align}
where $\boldsymbol{\mu}$ and $\boldsymbol{\Sigma}$ are defined as:
\begin{align*}
    \boldsymbol{\mu} = m(\boldsymbol{H^*}) + K(\boldsymbol{H^*}, \boldsymbol{H}) K(\boldsymbol{H}, \boldsymbol{H})^{-1}(\boldsymbol{E} - m (\boldsymbol{H})) \\ 
    \boldsymbol{\Sigma} = K(\boldsymbol{H^*}, \boldsymbol{H^*}) - K(\boldsymbol{H^*, \boldsymbol{H}})K(\boldsymbol{H}, \boldsymbol{H})^{-1}K(\boldsymbol{H}, \boldsymbol{H^*})
\end{align*}
$m(X)$ is the mean function, which represents prior knowledge about the function. In this, we encode the information about the previous $m$ hours. 
For our purpose, we utilize the mean of the previous days' energy prices to predict the mean of the next days' energy prices. For a time step, $k > i$, the prediction would look like,
\begin{align*}
    m(X_k) = \frac{\sum_{j = i-m}^{i}E_j[h_k == h_j]}{\sum_{j = i-m}^{i}1[h_k == h_j]}
\end{align*}

$K(X, X^*)$ is the kernel function, which controls smoothness of the function. It is constructed to measure the relation between different input values. Since energy usage is periodic in nature, with a period occurring every 24 hours, consideration must be taken for both the relationship of the periodic data while placing higher priority on recent data. With this in mind, we utilize the locally periodic kernel function from \cite{duvenaud_automatic_nodate}.
\begin{equation*}
    SQ_{exp}(x,x') = exp \left(-\frac{(x-x')^2}{2 \times l_{exp}^2}\right) 
\end{equation*}
\begin{equation}
    K(x, x') = \sigma^2exp \left(\frac{-2sin^2(\pi|x-x'|/\rho)}{l_{per}^2}\right)SQ_{exp}(x,x')
\end{equation}
where $l_{exp} = 24$, $l_{per} = \frac{3}{7}$, $\rho = 24$, $\sigma = 1000$. The $\rho$ value of $24$ illuminates that energy cycle usage have periods of 24 hours, and the $\sigma$ value of $1000$ accounts for the magnitude of values in the energy usage dataset we utilize which tracks energy usage in MW for every hour \cite{noauthor_hourly_nodate}. After fitting a Gaussian model to predict the next day's energy usage, we scale this proportionately to create a model for energy prices, $p$, based on simple supply-and-demand: $\hat{p}[\frac{\$}{kWh}] = \frac{\hat{E}}{80000}$ and $p[\frac{\$}{kWh}] = \frac{E}{80000}$. 

\begin{figure} [htbp]
    \centering
    \includegraphics[width=0.45\textwidth]{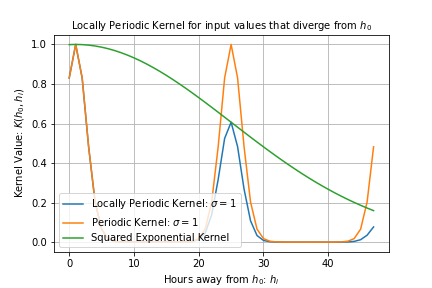}
    \caption{Demonstrates the general trend of the locally periodic kernel.}
    \label{fig:LocallyPeriodicKernel}
\end{figure}

This formulation allows for two important results:
\begin{enumerate}
    \item Uncertainty grows smallest when closest to the last recorded data (i.e. the prediction at i will be most confident).
    \item Uncertainty and mean have a periodic relation with values from the previous days
\end{enumerate}

\subsection{Curve Sampling}

Sampling is done according to Equation \ref{eq:post} with $\mathcal{N}(\boldsymbol{\mu}, \boldsymbol{\Sigma})$. In order to require a smoother fit, numpy's polynomial fitting toolbox \cite{harris2020array} is used to fit the sampled data to a $n^{th}$ order polynomial. This was attempted along with utilizing scipy's curve fitting toolbox \cite{2020SciPy-NMeth} to a sinusoidal function. Ultimately, polynomial fitting was chosen to allow for locally periodic behavior as modelled by the kernel definition. 

\begin{figure} [htbp]
    \centering
    \includegraphics[width=0.45\textwidth]{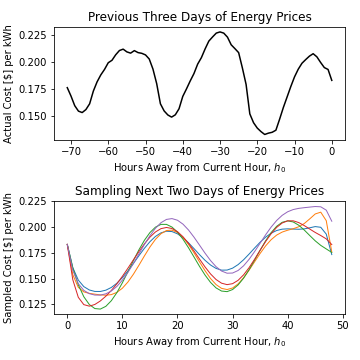}
    \caption{Given the previous three days of energy prices (depicted above), a Gaussian process is constructed using the defined Kernel: $K(X,X')$, and five samples were collected and plotted (below) according to $\mathcal{N}(\mu, \boldsymbol{\Sigma})$. Note the higher variance the farther from $h_0$ and the general periodic trend.}
    \label{fig:Gaussian_Process_Sampling}
\end{figure}

As can be seen in Figure \ref{fig:Gaussian_Process_Sampling}, samples for the next two days can deviate slightly but follow a locally periodic trend. The most variation occurs near $h_i = 48$, where the fewest measured values are near. Given the definition for the Kernel function, this is to be expected.

\subsection{Consumer Modelling}

\begin{figure} [htbp]
    \centering
    \includegraphics[width=0.4\textwidth]{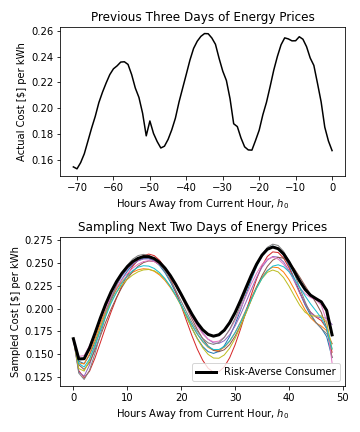}
    \caption{Risk-Averse Consumer Sampling chooses to optimize over the sample for the worst case situation: $n_s$ = 10}
    \label{fig:Gaussian_Process_Sampling3}
\end{figure}

There are two types of consumers to consider, the average consumer, and the risk-averse consumer. The average consumer will optimize over the predicted mean, $\boldsymbol{\mu}$, calculated from the Gaussian process while the risk-averse consumer will take $n_{s}$ samples and optimize over the worst case. An example of the risk-averse consumer is shown in Figure \ref{fig:Gaussian_Process_Sampling3}. 

\subsection{Optimization Methods}

The multi objective function $f(\textbf{P})$ as defined in Equation \ref{eq:objective-fn} is a non-linear and non-differentiable function. These properties largely arise because of the quadratic function to calculate $d_t$ and the experimentally determined piece-wise functions for $d_t^{1C}$ and $\psi_t$. Furthermore, absolute value functions create cusps at which the function is non-differentiable. As a result, optimization methods requiring analytical gradients cannot be used. While linear programming, which in many cases can find true global optima, could be viable if the model were linearized like in \cite{maheshwari_optimizing_2020}, doing so would obfuscate the nonlinear intricacies introduced by the model.

Moreover, the optimization method must consider the constraints on $P_t^b$ and $SOC_t$ discussed in the previous sections.

The model ultimately achieved best results using central difference Nesterov momentum gradient descent with a quadratic constraint penalty. While computationally expensive considering the high dimensional design space, numerical gradient calculation is quite accurate and relatively straightforward to apply even to non-linear and non-differentiable functions. The quadratic penalty method further simplifies the problem by transforming a constrained optimization problem to an unconstrained one. The resulting objective function is as follows:

\begin{equation}
    \min_{\textbf{P}} f(\textbf{P}) + \rho \times p_{q}(\textbf{P})
    \end{equation}

where $\textbf{P}$ is the $12$ dimensional design variable, $f(\textbf{P})$ is the original multi-objective as formulated in Equation \ref{eq:objective-fn}, $\rho$ is the penalty scaling factor and $p_{q}(\textbf{P})$ is the quadratic penalty calculation, which is formulated below.

\begin{dmath}
p_{q}(\textbf{P}) = \sum_{t}^{T} \max(SOC_t - 100, 0)^2 + \min(SOC_t, 0)^2 
+ \max(\lvert P_t^b \rvert - P^{b, max}, 0)^2
\end{dmath}

Through experimentation with hyperparameters, it was found that this optimization method converged in about 1000 iterations. To avoid becoming stuck at local optima, several trials were conducted with randomly initialized starting points.
\section{Results}

\subsection{Baseline}

\begin{figure} [htbp]
    \centering
    \includegraphics[width=0.45\textwidth]{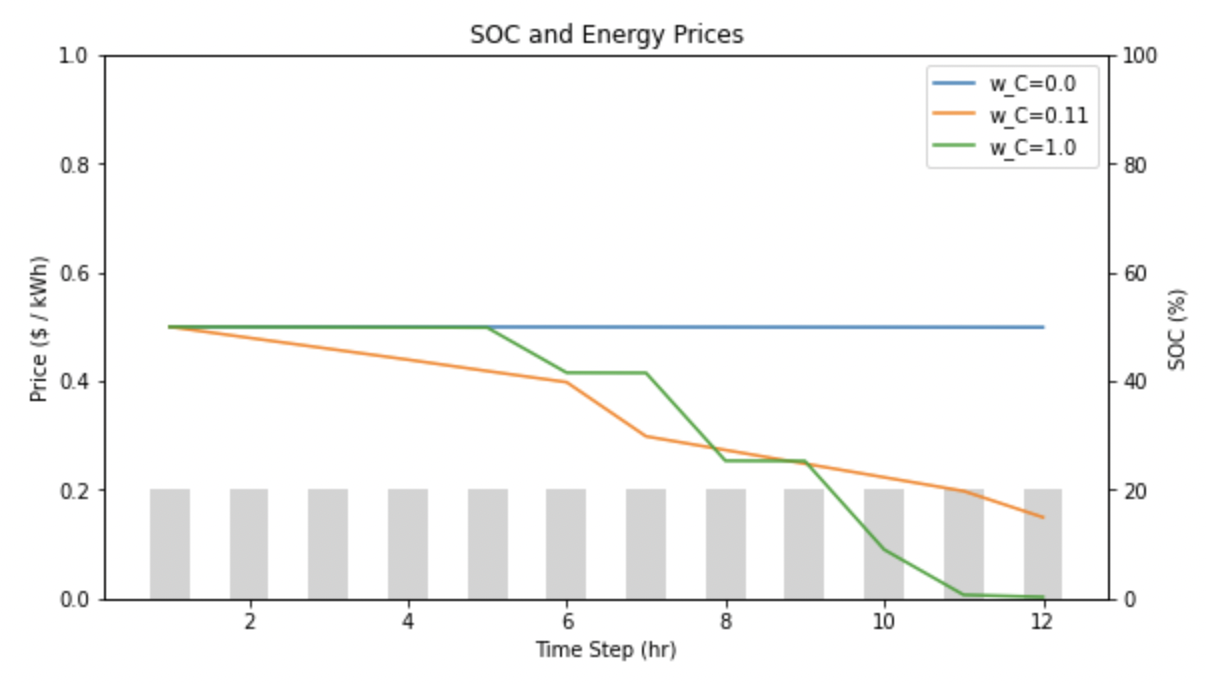}
    \caption{Baseline SOC schedule and energy prices with constant prices for three weighting factors. $w_C$ represents the weight placed on cost.}
    \label{fig:SOCandPricesConstPrices}
\end{figure}

The baseline method to which results will be compared is conducted assuming constant energy prices. The preliminary results are also a good indication of a properly implemented model and optimization method. These results are plotted in Figure \ref{fig:SOCandPricesConstPrices}.

Depending on the weighting factor $w_C$, the output charging schedule varies greatly. When $w_C=0.0$, the model reduces to a single objective optimization problem, with the only aim of minimizing battery degradation. As a result, the intuitive yet trivial proposed charging schedule is to never charge or discharge, resulting in zero degradation and a static SOC.

The other extreme, when $w_C=1.0$, transforms the model to a sole cost minimizing function. The proposed charging schedule results in completely discharging the battery by the end of the twelfth time interval. In order to incentivize energy usage in the model, negative cost is incurred (equivalent to gaining revenue) when the battery is discharged. The amount received is $\eta \times \lambda_t$ where $\eta$ is less than $1.0$ and $\lambda_t$ is the price of energy ($\frac{\$}{kWh}$) at time $t$. As such, under the assumption of constant prices, discharging and charging back to the same SOC will result in a positive cost. Thus, to gain the most revenue, the battery must completely discharge any existing energy and never charge during the 12 time steps. This behavior is represented in Figure \ref{fig:SOCandPricesConstPrices}.

Lastly, setting $w_C$ to an intermediate value between $0.0$ and $1.0$ results in a charging schedule that attempts to optimize for both degradation and cost.

It is important to note that multiple global optima can exist. For instance, once again assuming constant prices, if $w_C$ is set to $1.0$ then a charging schedule that discharges completely in one time step results in an equivalent cost outcome to a charging schedule that discharges one twelfth of the initial SOC at each of the twelve time steps.

\begin{figure} [htbp]
    \centering
    \includegraphics[width=0.45\textwidth]{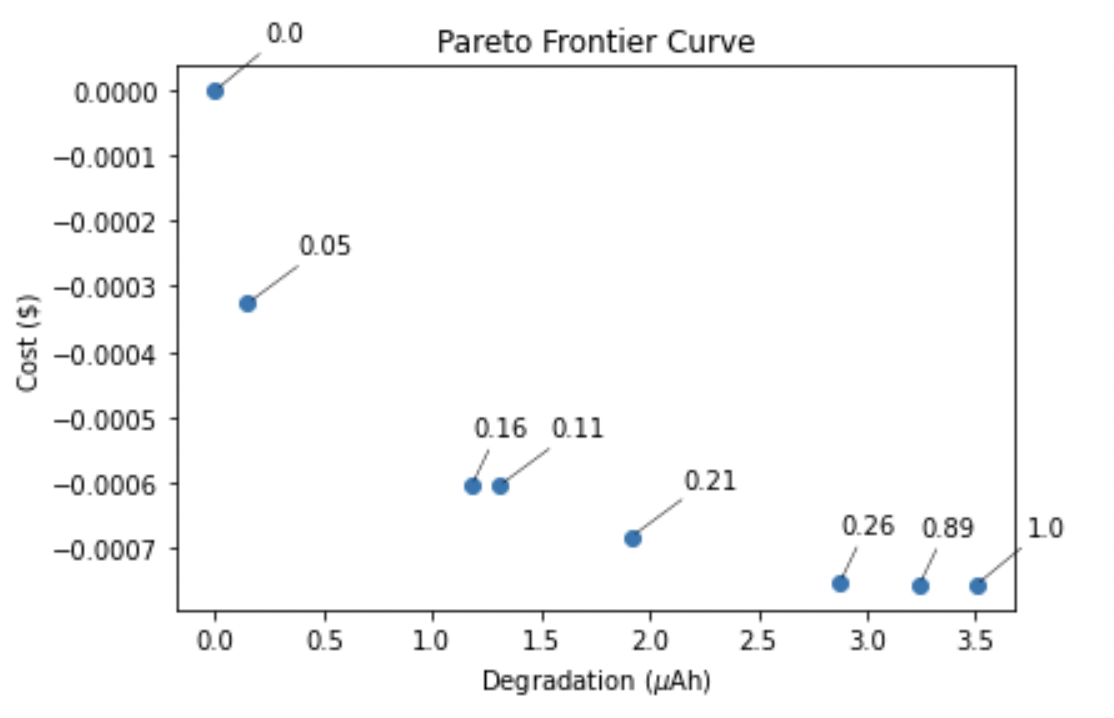}
    \caption{Baseline Pareto frontier with constant prices for various weighting factors. $w_C$ represents the weight placed on cost.}
    \label{fig:ParetoConstPrice}
\end{figure}

Apart from visualizing individual charging schedules proposed by the model, a Pareto frontier curve can be generated by sweeping through various weighting values corresponding to cost $C$ and degradation $D$. Figure \ref{fig:ParetoConstPrice} shows the curve corresponding to the optimal objective function $f(\textbf{P})$ evaluations for a given weight $w_C$ (and correspondingly a given weight $w_D=1.0-w_C$). The point corresponding to $w_C=0.0$ shows zero incurred cost and zero degradation, matching the respective constant charging schedule in Figure \ref{fig:SOCandPricesConstPrices}.

To formulate using economics terminology, points on the frontier lie on a curve of equal utility. That is, no points are dominated by any other point as defined by Kochenderfer \cite{mykel}. Moving along the curve is a matter of preference.

\subsection{Varied Prices}

While evaluating performance on a baseline model provides insight into how the model operates, it is not the most representative of true energy prices. The subsequent results relax the assumption of constant prices, producing proposed charging schedules and a Pareto frontier.

\begin{figure} [htbp]
    \centering
    \includegraphics[width=0.45\textwidth]{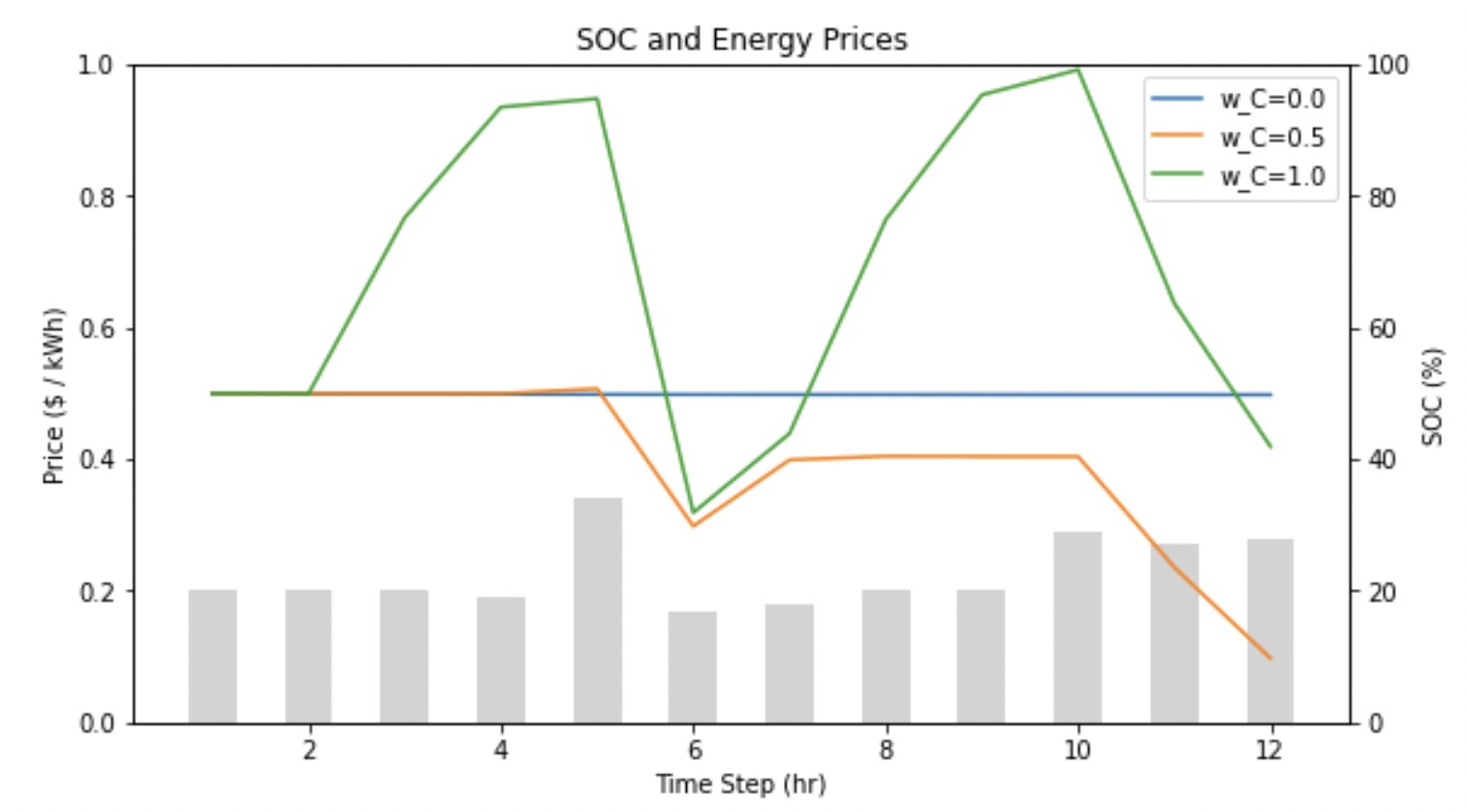}
    \caption{Pareto frontier with randomly varied prices for various weighting factors. $w_C$ represents the weight placed on cost.}
    \label{fig:SOCandPricesVarPrices}
\end{figure}

Figure \ref{fig:SOCandPricesVarPrices}, like with the baseline model, shows three proposed charging schedules for weights of $w_C = 0.0, 0.5,$ and $ 1.0$. For the $w_C=0.0$ case in which only degradation is optimized for, the model dictates to never charge or discharge. Regardless of what the energy prices may be, this will always be the optimal charging strategy with this weighting factor.

The $w_C=1.0$ case, in which only cost is optimized for, demonstrates extreme fluctuations in charge. More interestingly, the model takes advantage of the lower energy prices and suggests to maximally charge during these intervals. When energy prices increase, as in $t=5$, the battery is greatly discharged, likely limited only by maximum power output $P_{b,max}$. During intervals of high prices, discharging results in relatively higher revenue gained than when prices are low. Then when prices are low, the model encourages charging as a smaller cost is incurred.

The intermediate case of $w_C=0.5$ weights both cost and degradation equally. As a result, the proposed charging strategy is a less aggressive version than the one when $w_C=1.0$ and a more aggressive version than the one when $w_C=0.0$.

\begin{figure} [htbp]
    \centering
    \includegraphics[width=0.45\textwidth]{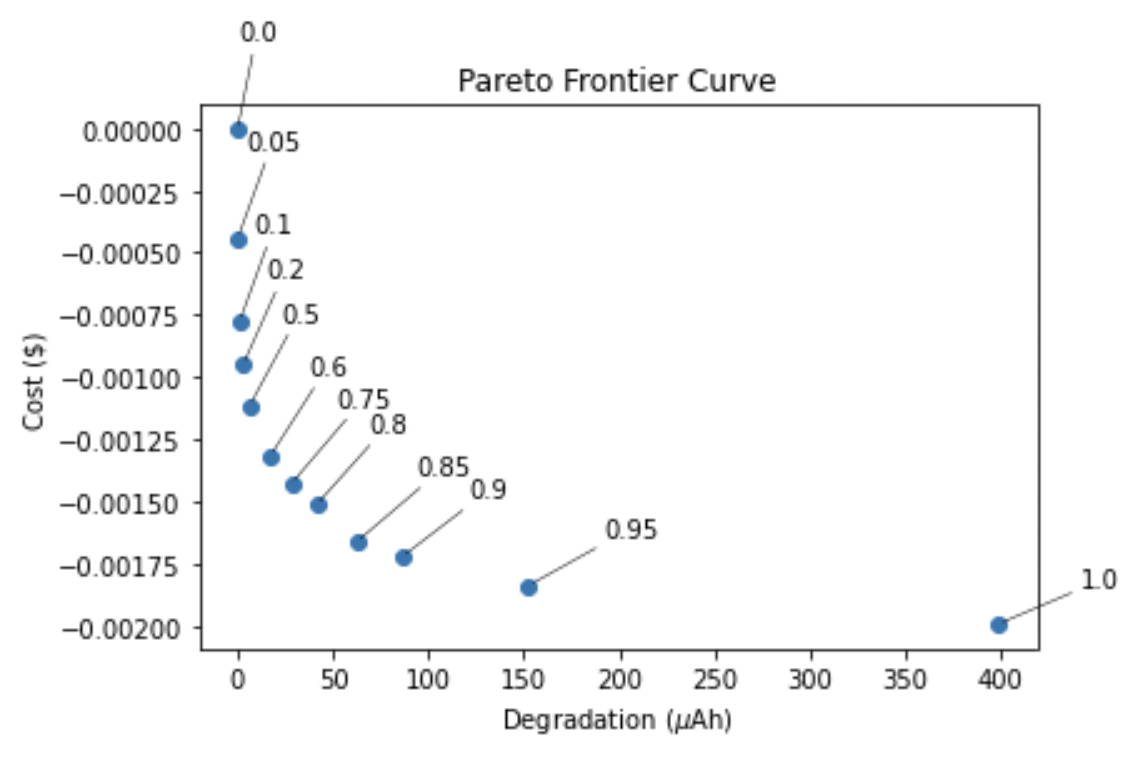}
    \caption{Pareto frontier with randomly varied prices for various weighting factors. $w_C$ represents the weight placed on cost.}
    \label{fig:ParetoVarPrice}
\end{figure}

Finally, the corresponding Pareto frontier can be seen in Figure \ref{fig:ParetoVarPrice}. 

\subsection{Predicted Prices}
\begin{figure} [htbp]
    \centering
    \includegraphics[width=0.45\textwidth]{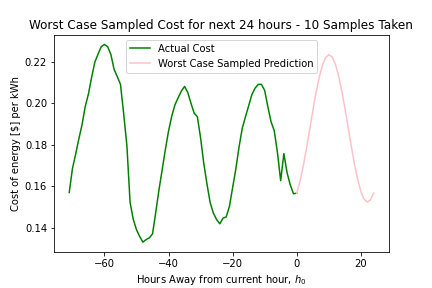}
    \caption{Worst Case sample for next 24 hours given previous 3 days of energy, $n_s = 10$}
    \label{fig:WorstCaseSample}
\end{figure}

Utilizing the Gaussian Process model and treating the consumer as a risk-averse consumer, we take $n_s = 10$ samples and optimize over the worst case scenario as shown in Figure \ref{fig:Gaussian_Process_Sampling3}. For this case, we optimize over the next 24 hours of predicted energy prices to emulate the real-world situation of having to optimize and plan for the full next day of energy. The sampled prediction is shown in Figure \ref{fig:WorstCaseSample}. Note, hour $h_0$ does not correlate to midnight. It is a randomly chosen hour on any given day.

\begin{figure} [htbp]
    \centering
    \includegraphics[width=0.45\textwidth]{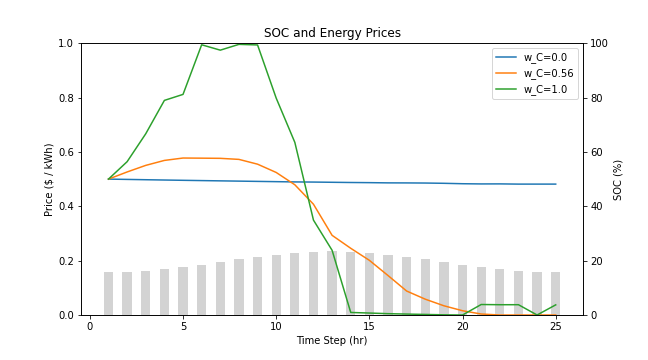}
    \caption{Optimization of Worst Case Sample, $n_s = 10$. Charging is shown to be promoted for higher $w_c$ values when future prices are going to be higher.}
    \label{fig:WorstCaseSOC}
\end{figure}

The results shown in Figure \ref{fig:WorstCaseSOC} demonstrate the extreme cases of $w_c = 1.0$ and $w_c = 0.0$, and it exemplifies the middle case of $w_c = 0.56$. Similar trends to the varied price schedule is demonstrated in Figure \ref{fig:SOCandPricesVarPrices}, where charging is considered when future prices will be higher. 

\begin{figure} [htbp]
    \centering
    \includegraphics[width=0.45\textwidth]{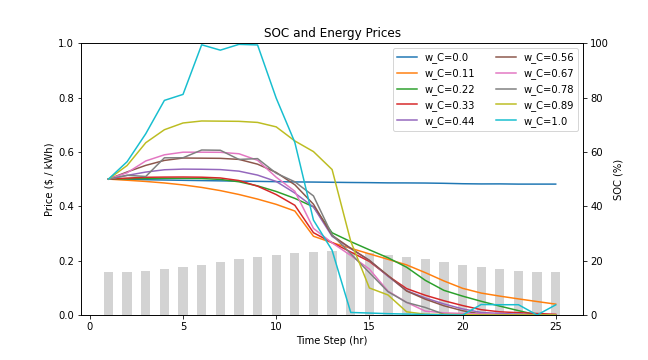}
    \caption{Optimization of Worst Case Sample for several $w_c$ values. Higher values of $w_c$ are shown to prioritize charging more aggressively while lower values are more cautious about battery degradation which matches the optimization function defined in Equation \ref{eq:objective-fn}.}
    \label{fig:WorstCaseSOCall}
\end{figure}

\begin{figure} [htbp]
    \centering
    \includegraphics[width=0.45\textwidth]{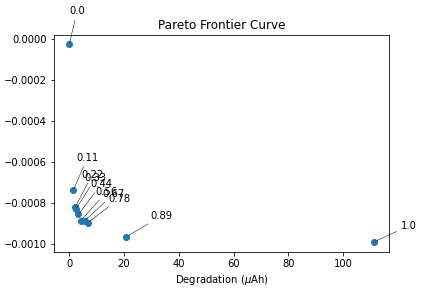}
    \caption{Pareto Frontier Curve for Worst Case Optimization over a 24-hour Trading Period.}
    \label{fig:WorstCasePareto}
\end{figure}

To exemplify the effect of changing $w_c$ even more drastically, Figure \ref{fig:WorstCaseSOCall} is provided. The Pareto frontier curve is also provided in Figure \ref{fig:WorstCasePareto}.

\section{Discussion}
From the results, it is evident that utilizing central difference Nesterov momentum gradient descent is applicable for a variety of different cost forecasts. For the situation in which costs remain constant, charging is inefficient if the revenue function is considered as "buying" and "selling back" to the market. Due to the $\eta$ inefficiency factor, $P_t^{b,ch} = \eta P_t^{m,ch}$ and  $\eta P_t^{b,dis} = P_t^{m,dis}$. Therefore, with constant prices, it never makes sense to charge, regardless of the $w_c$ value as revenue will be negative and the battery will degrade. Figure \ref{fig:SOCandPricesConstPrices} demonstrates this phenomenon, with the battery only discharging over the entire time interval.

For the varied price forecast and the worst-case realistic price forecast, charging is prioritized when future prices are predicted to be higher. An intuitive manner of understanding this is that the buyer is "buying" the energy at a lower price and "selling" it back to the market at a higher price. If the price difference is large enough to overcome the $\eta$ inefficiency factor, the buyer will still be encouraged to buy and then sell. This can be seen in Figure \ref{fig:SOCandPricesVarPrices} and Figure \ref{fig:WorstCaseSOC}. 

The Gaussian price modelling provides an important extension to the work of \cite{maheshwari_optimizing_2020}. By introducing uncertainty into the model, optimization can be done for battery charging in the future. The locally periodic kernel of \cite{duvenaud_automatic_nodate} proves to be an effective way to model the periodic, varying pattern of costs over 24-hour cycles. The sampling performed over Gaussian distributions for future costs simulates different cost forecasts, and different types of consumers can choose how to optimize given those samples.

\section{Conclusion}
This work extends the work done by \cite{maheshwari_optimizing_2020}. The degradation model used by \cite{maheshwari_optimizing_2020} in conjunction with the revenue model (with associated inefficiencies) is an effective manner to balance battery upkeep and encourage battery usage at optimal intervals. 

Instead of utilizing linearization techniques as \cite{maheshwari_optimizing_2020} does, this paper utilizes central difference Nesterov momentum gradient descent to come to optimal charging strategies. Despite the large state space and high computation time, this strategy proves effective in coming up with charging strategies for the next $n$ trading intervals given knowledge of the cost forecast. Utilizing baseline and randomly varied prices, charging strategies were demonstrated to match intuition. 

The cost forecast model was completed utilizing Gaussian Processes. The utilization of a locally periodic kernel captures the periodic, yet varying, structure of energy prices over 24-hour cycles. The formulated Gaussian Process distribution is then sampled over and fit to a polynomial in order to sample over the next $n$ trading intervals. This is shown to be an effective manner of capturing future uncertainty for prices. Different types of consumers can utilize the uncertainty of this specified Gaussian Process distribution to decide how to proceed.

The results of the SOC schedule and energy price visualizations and the Pareto frontier plots demonstrate the balance between battery degradation and the revenue model. Point along the Pareto frontier can be thought of as values of equal utility, in that no points are dominated by any other point. This specification further allows adaptability for persons who have varying priority of revenue vs. battery degradation depending on the use case. 

\section{Future Work}
 As mentioned in \cite{maheshwari_optimizing_2020}, more work can be done to create a battery degradation model as optimization focuses on battery usage have seldom considered this important factor of the equation. 
 
 Furthermore, this work can be extended upon to integrate further optimization over the generated Gaussian Process samples, instead of choosing just one sample and optimizing over that. Applying a gradient measurement strategy from \cite{mykel} could allow for optimization while directly measuring and adapting to the uncertainties of the future.

\section{Group Member Contributions} 
Jacob:
\begin{itemize}
    \item Implemented battery degradation model in Python
    \item Added functionality for quadratic penalty method
    \item Programmed finite differences calculation and gradient descent functions
    \item Tuned hyperparameters for optimal performance
    \item Produced graphs

\end{itemize}

\noindent Nico:
\begin{itemize}
    \item Experimented with kernels for Gaussian Process Model
    \item Implemented GPM for energy prices 
    \item Integrated GPM model with optimization methods
    \item Produced Graphs

\end{itemize}

\section{Appendix}

A link to our GitHub can be found  \href{https://github.com/jacobazoulay/AA222_Final_Project}{here}.

\bibliographystyle{ieeetr}
\bibliography{references.bib}

\end{document}